\documentclass[11pt,reqno]{amsart}
\usepackage{amsmath, amssymb, amsthm}
\usepackage{url}
\usepackage[breaklinks]{hyperref}
\usepackage{csquotes}

\renewcommand{\(}{\left\(}
\renewcommand{\)}{\right\)}
\renewcommand{\[}{\left\[}
\renewcommand{\]}{\right\]}

\renewcommand{\mkbegdispquote}[2]{\itshape}
\numberwithin{equation}{section}
 \theoremstyle{plain}
\newtheorem{theorem}{Theorem}[section]

\newtheorem{remark}[]{Remark}

   \makeatletter
\def\proof{\@ifnextchar[{\@oproof}{\@nproof}}
\def\@oproof[#1][#2]{\trivlist\item[\hskip\labelsep\textit{#2 Proof of\
#1.}~]\ignorespaces}
\def\@nproof{\trivlist\item[\hskip\labelsep\textit{Proof.}~]\ignorespaces}

\makeatother

\begin{document}
\title[Infinite families of solutions for $A^3 + B^3 = C^3 + D^3$ and $A^4 + B^4 + C^4 + D^4 + E^4 = F^4$]{Infinite families of solutions for $A^3 + B^3 = C^3 + D^3$ and $A^4 + B^4 + C^4 + D^4 + E^4 = F^4$  in the spirit of Ramanujan}

\author{Archit Agarwal}
\address{Archit Agarwal, Department of Mathematics, Indian Institute of Technology Indore, Simrol, Indore 453552, Madhya Pradesh, India.}
\email{phd2001241002@iiti.ac.in, archit.agrw@gmail.com}

\author{Meghali Garg}
\address{Meghali Garg, Department of Mathematics, Indian Institute of Technology Indore, Simrol, Indore 453552, Madhya Pradesh, India.}
\email{phd2001241005@iiti.ac.in, meghaligarg.2216@gmail.com}

\thanks{2010\textit{Mathematics Subject Classification.} Primary 11D25.}

\maketitle
\begin{abstract}

Ramanujan,  in his lost notebook,  gave an interesting identity, which generates infinite families of solutions to Euler's Diophantine equation $A^3 + B^3 = C^3 + D^3$.  In this paper,  we produce a few infinite families of solutions to the aforementioned Diophantine equation as well as for the Diophantine equation $A^4 + B^4 + C^4 + D^4 + E^4 = F^4$ in the spirit of Ramanujan.


\end{abstract}

\section{Introduction}

Littlewood once remarked that \emph{``Every natural number was one of Ramanujan's personal friends''.} $1729$ is one such special number which was dear to Ramanujan. This number is popularly known as Ramanujan taxicab number or Hardy-Ramanujan number due to the following incident.

The story begins in 1918, when Ramanujan was admitted to Matlock Sanitorium in Derbyshire. Hardy went to see him and communicated the following story about his visit:
\begin{quote}
\emph{I remember once going to him when he was lying ill at Putney. I had ridden in taxi cab number 1729 and remarked that the number seemed to me rather a dull one, and that I hoped it was not an unfavourable omen. No, he replied, it is a very interesting number; it is the smallest number expressible as the sum of two cubes in two different ways.}
\end{quote}

One can easily verify that $1729$ is the smallest natural number to be represented as sum of two cubes in two different ways i.e., $1729=1^3+12^3=9^3+10^3$. Once we know this special property of $1729$, this prompts one to think are there more numbers satisfying this property? This is analogous to solving the Diophantine equation
\begin{align}\label{Euler's Diophantine equation}
A^3+B^3=C^3+D^3,
\end{align}
which is well-known as Euler's Diophantine equation. Ramanujan, in his lost notebook \cite[p.~341]{lost notebook}, mentioned the following amazing identity, which gives infinitely many solutions to \eqref{Euler's Diophantine equation}.
\begin{theorem}\label{ramanujan main result}
If
\begin{align}\label{ramanujan identity}
\begin{split}
\frac{1+53x+9x^2}{1-82x-82x^2+x^3}=\sum_{n=0}^\infty a_nx^n=\sum_{n=1}^\infty \alpha_{n-1}x^{-n},\\
\frac{2-26x-12x^2}{1-82x-82x^2+x^3}=\sum_{n=0}^\infty b_nx^n=\sum_{n=1}^\infty \beta_{n-1}x^{-n},\\
\frac{2+8x-10x^2}{1-82x-82x^2+x^3}=\sum_{n=0}^\infty c_nx^n=\sum_{n=1}^\infty \gamma_{n-1}x^{-n},
\end{split}
\end{align}
then\footnote{In \cite[p.~203, Equation (8.5.13)]{LN4}, the authors mentioned that the second equality of \eqref{ramanujan final result} is incorrect. However, we emphasize that the second equality is in fact correctly mentioned by Ramanujan.}
\begin{align}\label{ramanujan final result}
a_n^3+b_n^3=c_n^3+(-1)^{n} \hspace{5mm} \textrm{and} \hspace{5mm} \alpha_n^3+\beta_n^3=\gamma_n^3+(-1)^n.
\end{align}
\end{theorem}

A few solutions obtained by Ramanujan corresponding to $a_n^3+b_n^3=c_n^3+(-1)^{n}$, for $n=1,2$, are 
\begin{align*}
135^3 + 138^3 &= 172^3 - 1, \\ 11161^3 + 11468^3 &=  14258^3 + 1, 
\end{align*}
and from $\alpha_n^3+\beta_n^3=\gamma_n^3+(-1)^n$, for $n=0,1,2$, he derived
\begin{align*}
9^3 + (-12)^3 &= (-10)^3 +1, \\ 791^3 + (-1010)^3 &= (-812)^3 -1, \\ 65601^3 + (-83802)^3 &= (-67402)^3 +1.
\end{align*}
Like many of Ramanujan's results, one wonders how he managed to achieve such an identity.
 Hirschhorn examined Ramanujan's claim over a period of time and put forward the idea about how Ramanujan might have proceeded in \cite{hir 2}-\cite{hir 4}.  Hirschhorn along with Han \cite{hir 1} gave an alternative proof of Theorem \ref{ramanujan main result}. McLaughlin \cite{JMc} found an identity similar to \eqref{ramanujan identity} involving eleven sequences. Over the years, this magical number $1729$ has gained attention of many mathematicians. While trying to find more integers satisfying Euler's Diophantine equation, Silverman \cite{silverman} tackled the problem using elliptic curves. Chen \cite{chen} discussed an algorithm to find various identities similar to \eqref{ramanujan identity}. A few years back, Ni\c{t}ic\u{a} had observed a very interesting property about $1729$. It says that if we add the digits of $1729$, we get $19$, then multiplying $19$ with $91$, which is obtained by reversing the digits of $19$, we again get $1729$. The only integers satisfying this property are $\{1, 81, 1458, 1729\}$. In \cite{VN}, Ni\c{t}ic\u{a} obtained a generalization of above property. An elegant generalization of taxicab number has been discussed by Dinitz, Games and Roth \cite{DGR} in 2019, where they discussed the smallest number that can be represented as sum of $n$ positive $m^{th}$ powers in atleast $t$ ways.

In this paper, we use Hirschhorn's idea to obtain infinite families of solutions to \eqref{Euler's Diophantine equation}. Hirschhorn believes that Ramanujan might have used his previously obtained parametric solution of \eqref{Euler's Diophantine equation} to get \eqref{ramanujan identity}. One of Ramanujan's families of solutions \cite[p.~56]{RN4},  \cite[p.~266]{ramanujantifr}, \cite{ramanujanquestion} to \eqref{Euler's Diophantine equation} was given by 
\begin{align}\label{parametric_solution}
(3a^2+5ab-5b^2)^3+(4a^2-4ab+6b^2)^3+(5a^2-5ab-3b^2)^3=(6a^2-4ab+4b^2)^3.
\end{align}
On suggestion by Craig, Hirschhorn changed the variables as $a=A+B$ and $b=A-2B$ in $\eqref{parametric_solution}$ and obtained the following identity:
\begin{align}\label{after_variable_change}
(A^2+7AB-9B^2)^3+(2A^2-4AB+12B^2)^3=(2A^2+10B^2)^3+(A^2-9AB-B^2)^3.
\end{align}
We use \eqref{after_variable_change} to get various identities similar to \eqref{ramanujan identity} which leads to several families of solutions to \eqref{Euler's Diophantine equation}.
In the next section, we mention our main results. 
\section{Main Results}
At the beginning of Ramanujan's second notebook \cite[p.~3]{ramanujantifr}, he mentioned the following parametric solution to \eqref{Euler's Diophantine equation}.
\begin{theorem}\label{general parametric solution}
If $p^3+q^3+r^3=s^3$ and $m=  (s+q)\sqrt{\frac{s-q}{r+p}},~n= (r-p)\sqrt{\frac{r+p}{s-q}}$, then for any arbitrary $a,b$, we have
\begin{align}
(pa^2+mab-rb^2)^3+(qa^2-nab+sb^2)^3+(ra^2-mab-pb^2)^3=(sa^2-nab+qb^2)^3.
\end{align} 
\end{theorem}
Berndt \cite[p.~54]{RN4} gave a proof of the above theorem.  In this paper,  we present an elementary proof,  which also motivated us to  obtain the below identity.  

\begin{theorem}\label{6_variable_parametric_solution}
If $p^3+q^3+r^3+s^3+t^3=u^3$ and  $g=  (q-t)\sqrt{\frac{q+t}{p+s}},~h= (p-s)\sqrt{\frac{p+s}{r-u}},~k= (p-s)\sqrt{\frac{p+s}{q+t}},~l= (q-t)\sqrt{\frac{q+t}{r-u}},~m=  (r+u)\sqrt{\frac{r-u}{p+s}},~n=  (r+u)\sqrt{\frac{r-u}{q+t}}$, then for any arbitrary $a,b,c$, we have
\begin{align}
&(sa^2+pb^2+pc^2-gab-mac)^3 + (qa^2+tb^2+qc^2-kab-nbc)^3 + (ra^2+rb^2-uc^2-hac-lbc)^3 \nonumber  \\
 & +(pa^2+sb^2+sc^2+mac+gab)^3 + (ta^2+qb^2+tc^2+kab+nbc)^3 = (ua^2+ub^2-rc^2-hac-lbc)^3.
\end{align}
\end{theorem}
Next identities are in the spirit of Theorem \ref{ramanujan main result}.  
\begin{theorem}\label{fibonacci}
Suppose we have the following power series:
\begin{align}
\begin{split}
\frac{1-3x+9x^2}{1-2x-2x^2+x^3}=\sum_{n=0}^\infty a_nx^n,\\
\frac{2+6x-12x^2}{1-2x-2x^2+x^3}=\sum_{n=0}^\infty b_nx^n,\\
\frac{2+8x-10x^2}{1-2x-2x^2+x^3}=\sum_{n=0}^\infty c_nx^n,\\
\frac{1-11x+x^2}{1-2x-2x^2+x^3}=\sum_{n=0}^\infty d_nx^n,
\end{split}
\end{align}
 then, for all $n \geq 0$, we get
\begin{align}
a_n^3+b_n^3=c_n^3+d_n^3.
\end{align}
\end{theorem} 
A few solutions to the Euler's diophantine equation obtained from the above theorem for $n=1,2,3,4$ are:
\begin{align*}
(-1)^3+10^3&=12^3+(-9)^3,\\ 9^3+12^3&=18^3+(-15)^3,\\ 15^3+42^3&=58^3+(-49)^3,\\ 49^3+98^3&=140^3+(-119)^3
\end{align*}
\begin{remark} From the above examples, one can observe that $d_n=-a_{n+1}$ for all $n \geq 0$. We can easily verify this observation with the help of generating functions for $d_n$ and $a_n$.
\end{remark}

\begin{theorem}\label{gen -9^3n}
Let $a_n,~b_n,~c_n$ be the sequence of integers and $\alpha_n,~\beta_n,~\gamma_n$ be the sequence of rationals with the following generating functions:
\begin{align}\label{generating functions 9^n}
\begin{split}
\frac{2-8x-90x^2}{1-58x-522x^2+729x^3}=\sum_{n=0}^\infty a_nx^n=\sum_{n=1}^\infty \alpha_{n-1} x^{-n},\\
\frac{1+53x+9x^2}{1-58x-522x^2+729x^3}=\sum_{n=0}^\infty b_nx^n=\sum_{n=1}^\infty \beta_{n-1} x^{-n},\\
\frac{2+22x-108x^2}{1-58x-522x^2+729x^3}=\sum_{n=0}^\infty c_nx^n=\sum_{n=1}^\infty \gamma_{n-1} x^{-n},
\end{split}
\end{align} then one has 
\begin{align}\label{abc}
a_n^3+b_n^3=c_n^3+((-9)^n)^3 \hspace{5mm} \textrm{and} \hspace{5mm} \alpha_n^3+\beta_n^3=\gamma_n^3- \delta_n^3,
\end{align} where $\delta_n=(-1/9)^n$. 
\end{theorem}
\begin{remark}
We can get integral solutions from the second equality of \eqref{abc} on multiplying $\alpha_{n},~\beta_{n},~\gamma_{n}$ by $9^{2(n+1)}$. This will be clearer from the following example. For $n=1$, we have
\begin{align*}
\left(\frac{-10}{81}\right)^3 + \left(\frac{1}{81}\right)^3 &= \left(\frac{-4}{27}\right)^3 - \left(\frac{-1}{9}\right)^3.
\end{align*}
On clearing the denominators, we get $$(-10)^3+1^3=(-12)^3-(-9)^3.$$ Here we see that each of $\alpha_0,~\beta_0$ and $\gamma_0$ has been multiplied by $9^2$.
\end{remark}
A few more solutions obtained from the above theorem satisfying $a_n^3+b_n^3=c_n^3+((-9)^n)^3$, for $n=1,2,3$,  are
\begin{align*}
108^3 + 111^3 &= 138^3 + (-9)^3, \\ 7218^3 + 6969^3 &= 8940^3 + 81^3, \\ 473562^3 + 461415^3 &= 589098^3 + (-729)^3.
\end{align*}
Some more examples satisfying $\alpha_n^3+\beta_n^3=\gamma_n^3- \delta_n^3$ for $n=2,3$, after clearing the denominators, are
\begin{align*}
(-652)^3+535^3&=(-498)^3-81^3, \\ (-41578)^3+32281^3 &=(-33690)^3-(-729)^3.
\end{align*}

\begin{theorem}\label{gen 6^3n}
Suppose we have the following power series expansions:
\begin{align}
\begin{split}
\frac{2+22x+60x^2}{1+2x-12x^2-216x^3}=\sum_{n=0}^\infty a_nx^n=\sum_{n=1}^\infty \alpha_{n-1} x^{-n},\\
\frac{1-13x-6x^2}{1+2x-12x^2-216x^3}=\sum_{n=0}^\infty b_nx^n=\sum_{n=1}^\infty \beta_{n-1} x^{-n},\\
\frac{1+11x-54x^2}{1+2x-12x^2-216x^3}=\sum_{n=0}^\infty c_nx^n=\sum_{n=1}^\infty \gamma_{n-1} x^{-n},
\end{split}
\end{align} then we have
\begin{align}\label{2.6n}
a_n^3+b_n^3=c_n^3+(2\times 6^n)^3 \hspace{5mm} \textrm{and} \hspace{5mm} \alpha_n^3+\beta_n^3=\gamma_n^3- \delta_n^3,
\end{align} where $\delta_n=2/6^n$. 
\end{theorem}

\begin{remark}
Here again on multiplying by $6^{2(n+1)}$, rational solutions corresponding to $\alpha_{n},~\beta_{n},~\gamma_{n}$ can be converted to integral solutions of \eqref{Euler's Diophantine equation}. 
\end{remark}

A few more solutions to \eqref{Euler's Diophantine equation} extracted from Theorem \ref{gen 6^3n} are 
\begin{align*}
18^3 + (-15)^3 &= 9^3 + 12^3, \\ 48^3 + 36^3 &= (-60)^3 + 72^3, \\ 552^3 + (-36)^3 &= 444^3 + 216^3, \\ 3360^3 + (-2736)^3 &= 336^3 + 1296^3.
\end{align*}
The above solutions satisfy the first equality of \eqref{2.6n} and in accordance with the second equality, we clear the denominators to get
\begin{align*}
(-112)^3+76^3 &=(-84)^3-72^3, \\ (-328)^3+(-356)^3&=60^3-432^3.
\end{align*}


\begin{theorem}\label{Fibonacci_power_4_s_t}
Consider six different sequence of integers $\{a_n\}, \{b_n\}, \{c_n\}, \{d_n\}, \{e_n\},$ and $\{f_n\}$ with the following generating function:
\begin{align}\label{fibonacci a_n...f_n}
\begin{split}
\frac{8+8x+24x^2}{1-2x-2x^2+x^3}=\sum_{n=0}^\infty a_nx^n,\\
\frac{6-68x+18x^2}{1-2x-2x^2+x^3}=\sum_{n=0}^\infty b_nx^n,\\
\frac{14-60x+42x^2}{1-2x-2x^2+x^3}=\sum_{n=0}^\infty c_nx^n,\\
\frac{9+18x-27x^2}{1-2x-2x^2+x^3}=\sum_{n=0}^\infty d_nx^n,\\
\frac{4+8x-12x^2}{1-2x-2x^2+x^3}=\sum_{n=0}^\infty e_nx^n,\\
\frac{15+30x-45x^2}{1-2x-2x^2+x^3}=\sum_{n=0}^\infty f_nx^n.
\end{split}
\end{align} Then, we have $$a_n^4+b_n^4+c_n^4+d_n^4+e_n^4 = f_n^4.$$
\end{theorem}
A few numerical examples drawn from above theorem for $n=0,1,2$ are
\begin{align*}
8^4 + 6^4 + 14^4 + 9^4 + 4^4 &= 15^4,\\
24^4 + (-56)^4 + (-32)^4 + 36^4 + 16^4 &= 60^4,\\
88^4 + (-82)^4 + 6^4 + 63^4 + 28^4 &= 105^4,\\
\end{align*}

\begin{theorem}\label{3^n.8_power_4}
Consider the following power series
\begin{align}
\begin{split}
\frac{6+184x+54x^2}{1-28x-84x^2+27x^3}=\sum_{n=0}^\infty a_nx^n,\\
\frac{14-64x+126x^2}{1-28x-84x^2+27x^3}=\sum_{n=0}^\infty b_nx^n,\\
\frac{9-81x^2}{1-28x-84x^2+27x^3}=\sum_{n=0}^\infty c_nx^n,\\
\frac{4-36x^2}{1-28x-84x^2+27x^3}=\sum_{n=0}^\infty d_nx^n,\\
\frac{15-135x^2}{1-28x-84x^2+27x^3}=\sum_{n=0}^\infty e_nx^n,
\end{split}
\end{align}
then the coefficients of the above power series satisfy
\begin{align}
a_n^4+b_n^4+c_n^4+d_n^4 = e_n^4 - (8 \times (-3)^n)^4.
\end{align}
\end{theorem}
A few solutions obtained from above theorem corresponding to $n = 1,2,3$ are
\begin{align*}
352^4 + 328^4 + 252^4 + 112^4 &= 420^4 - (-24)^4,\\
10414^4 + 10486^4 + 7731^4 + 3436^4 &= 12885^4 - 72^4,\\
320998^4 + 320782^4 + 237393^4 + 105508^4 &= 395655^4 - (-216)^4,\\
\end{align*}

\begin{theorem}\label{fibonacci_power_4_m_n}
Let $\{a_n\}, \{b_n\}, \{c_n\}, \{d_n\}, \{e_n\},$ and $\{f_n\}$ be the sequences of integers with the following power series:
\begin{align}
\begin{split}
\frac{4-16x+12x^2}{1-2x-2x^2+x^3}=\sum_{n=0}^\infty a_nx^n,\\
\frac{3+6x-9x^2}{1-2x-2x^2+x^3}=\sum_{n=0}^\infty b_nx^n,\\
\frac{2-20x+6x^2}{1-2x-2x^2+x^3}=\sum_{n=0}^\infty c_nx^n,\\
\frac{4+8x-12x^2}{1-2x-2x^2+x^3}=\sum_{n=0}^\infty d_nx^n,\\
\frac{2+4x+6x^2}{1-2x-2x^2+x^3}=\sum_{n=0}^\infty e_nx^n,\\
\frac{5+10x-15x^2}{1-2x-2x^2+x^3}=\sum_{n=0}^\infty f_nx^n.
\end{split}
\end{align} Then, $$a_n^4+b_n^4+c_n^4+d_n^4+e_n^4 = f_n^4$$
\end{theorem}
A few solutions obtained from  the above theorem for $n=1,2,3$ are
\begin{align*}
(-8)^4 + 12^4 + (-16)^4 + 16^4 + 8^4 &= 20^4,\\
4^4 + 21^4 + (-22)^4 + 28^4 + 26^4 &= 35^4,\\
(-12)^4 + 63^4 + (-78)^4 + 84^4 + 66^4 &= 105^4.
\end{align*}

\begin{theorem}\label{2.3^n_Power_4}
Consider the following power series 
\begin{align}
\begin{split}
\frac{4-24x+36x^2}{1-39x-117x^2+27x^3}=\sum_{n=0}^\infty a_nx^n,\\
\frac{3-27x^2}{1-39x-117x^2+27x^3}=\sum_{n=0}^\infty b_nx^n,\\
\frac{4-36x^2}{1-39x-117x^2+27x^3}=\sum_{n=0}^\infty c_nx^n,\\
\frac{2+60x+18x^2}{1-39x-117x^2+27x^3}=\sum_{n=0}^\infty d_nx^n,\\
\frac{5-45x^2}{1-39x-117x^2+27x^3}=\sum_{n=0}^\infty e_nx^n,\\
\end{split}
\end{align}
then we have
\begin{align}
a_n^4+b_n^4+c_n^4+d_n^4 = e_n^4 - (2 \times (-3)^n)^4.
\end{align}
\end{theorem}
A few examples obtained from the above theorem are 
\begin{align*}
\begin{split}
132^4 + 117^4 + 156^4 + 138^4 &= 195^4 - (-6)^4,\\
5652^4 + 4887^4 + 6516^4 + 5634^4 &= 8145^4 - 18^4,\\
235764^4 + 204201^4 + 272268^4 + 235818^4 &= 340335^4 - (-54)^4.
\end{split}
\end{align*}
From these numerical evidences, one can clearly see that the families of solutions of each theorem are different.  In the next section, we present the proofs of all the results.
\section{Proof of Main Results}
\begin{proof}[Theorem \ref{general parametric solution}][]
Let $(p,q,r,s)$ be a known solution to $A^3+B^3+C^3=D^3$. Then consider a straight line passing through $(p,q,r,s)$ as 
\begin{align*}
\frac{A-p}{a}=\frac{B-q}{b}=\frac{C-r}{-a}=\frac{D-s}{b}=\theta.
\end{align*}
This gives
\begin{equation}\label{A,B,C,D}
\begin{aligned}
A&=a\theta+p,\\
B&=b\theta+q,\\
C&=-a\theta+r,\\
D&=b\theta+s.
\end{aligned}
\end{equation}
Substituting \eqref{A,B,C,D} in $A^3+B^3+C^3=D^3$, we get
\begin{align}\label{value of theta}
\theta=-\frac{a(p^2-r^2)+b(q^2-s^2)}{a^2(p+r)+b^2(q-s)}.
\end{align}
Plugging the value \eqref{value of theta} of $\theta$ in \eqref{A,B,C,D}, we get
\begin{align*}
A&=a^2r(p+r)-ab(q^2-s^2)-b^2p(s-q),\\
B&=a^2q(p+r)-ab(p^2-r^2)+b^2s(s-q),\\
C&=a^2p(p+r)+ab(q^2-s^2)-b^2r(s-q),\\
D&=a^2s(p+r)-ab(p^2-r^2)+b^2q(s-q).
\end{align*}
This is the most general solution. Now, we make a change of variable
\begin{align*}
a\rightarrow \frac{a}{\sqrt{p+r}},~~~b\rightarrow\frac{b}{\sqrt{s-q}}.
\end{align*}
Without loss of generality, we may assume $s>q$. This gives us
\begin{align*}
A&=a^2r+ab(s+q)\sqrt{\frac{s-q}{p+r}}-b^2p,\\
B&=a^2q-ab(p-r)\sqrt{\frac{p+r}{s-q}}+b^2s,\\
C&=a^2p-ab(s+q)\sqrt{\frac{s-q}{p+r}}-b^2r,\\
D&=a^2s-ab(p-r)\sqrt{\frac{p+r}{s-q}}+b^2q.
\end{align*}
Interchanging $r$ and $p$ gives Ramanujan's solution.
\end{proof}
\begin{proof}[Theorem $\ref{6_variable_parametric_solution}$][]
Let us consider a straight line passing through a point $(p,q,r,s,t,u)$, which is a known solution to $A^3+B^3+C^3+D^3+E^3 = F^3$,  
\begin{align*}
\frac{A-p}{a}=\frac{B-q}{b}=\frac{C-r}{c}=\frac{D-s}{-a}=\frac{E-t}{-b}=\frac{F-u}{c}=\theta.
\end{align*}
This gives
\begin{align}\label{all x_i's}
\begin{split}
A&=a\theta + p,  \\
B&=b\theta + q,  \\
C&=c\theta + r,  \\
D&=-a\theta + s,  \\
E&=-b\theta + t,  \\
F&=a\theta + u.
\end{split}
\end{align}
Substitute \eqref{all x_i's} in  $A^3+B^3+C^3+D^3+E^3 = F^3$ to get
\begin{align*}
\theta = \frac{-[a(p^2-s^2)+b(q^2-t^2)+c(r^2-u^2)]}{a^2(p+s)+b^2(q+t)+c^2(r-u)}.
\end{align*}
Utilizing the value of $\theta$ in \eqref{all x_i's}, we get
\begin{align*}
A&=a^2s(p+s)-ab(q^2-t^2)-ac(r^2-u^2)+b^2p(q+t)+c^2p(r-u),\\
B&=-ab(p^2-s^2)+b^2t(q+t)-bc(r^2-u^2)+a^2q(p+s)c^2q(r-u),\\
C&=-ac(p^2-s^2)-bc(q^2-t^2)-c^2u(r-u)+a^2r(p+s)+b^2r(q+t),\\
D&=a^2p(p+s)+ab(q^2-t^2)+ac(r^2-u^2)+b^2s(q+t)+c^2s(r-u),\\
E&=ab(p^2-s^2)+b^2(q+t)+bc(r^2-u^2)+a^2t(p+s)+c^2t(r-u),\\
F&=-ac(p^2-s^2)-bc(q^2-t^2)-c^2r(r-u)+a^2u(p+s)+b^2u(q+t).
\end{align*}
Now we make a change of variable
\begin{align*}
a\rightarrow \frac{a}{\sqrt{p+s}}, \hspace{3mm} b\rightarrow \frac{b}{\sqrt{q+t}}, \hspace{3mm} c\rightarrow \frac{c}{\sqrt{r-u}}.
\end{align*}
Without loss of generality, we may assume $r>u$. This gives us
\begin{align*}
A&=a^2s+b^2p+c^2p-ab(q-t)\sqrt{\frac{q+t}{p+s}}-ac(r+u)\sqrt{\frac{r-u}{p+s}},\\
B&=a^2q+b^2t+c^2q-ab(p-s)\sqrt{\frac{p+s}{q+t}}-bc(r+u)\sqrt{\frac{r-u}{q+t}},\\
C&=a^2r+b^2r-c^2u-ac(p-s)\sqrt{\frac{p+s}{r-u}}-bc(q-t)\sqrt{\frac{q+t}{r-u}},\\
D&=a^2p+b^2s+c^2s+ac(r+u)\sqrt{\frac{r-u}{p+s}}+ab(q-t)\sqrt{\frac{q+t}{p+s}},\\
E&=a^2t+b^2q+c^2t+ab(p-s)\sqrt{\frac{p+s}{q+t}}+bc(r+u)\sqrt{\frac{r-u}{q+t}},\\
F&=a^2u+b^2u-c^2r-ac(p-s)\sqrt{\frac{p+s}{r-u}}-bc(q-t)\sqrt{\frac{q+t}{r-u}}.
\end{align*}
This completes the proof.
\end{proof}
\begin{proof}[Theorem $\ref{fibonacci}$][]
Let us assume our sequence $\omega_n$ to be Fibonacci sequence defined as
\begin{align}\label{Fibonnaci_sequence}
\omega_{n+2}=\omega_{n+1}+\omega_n\hspace{3mm}\textrm{with}\hspace{3mm} \omega_0=0\hspace{3mm} \textrm{and} \hspace{3mm} \omega_1=1.
\end{align}
Then for the coexistence of $\omega_n$ with \eqref{after_variable_change}, we substitute
\begin{align*}
A=\omega_{n+1}\hspace{3mm} \textrm{and} \hspace{3mm}B=\omega_n,
\end{align*}
with
\begin{align}\label{a_n, b_n, c_n, d_n}
\begin{split}
&a_n=A^2+7AB-9B^2=\omega_{n+1}^2+7\omega_n\omega_{n+1}-9\omega_n^2,\\
&b_n=2A^2-4AB+12B^2=2\omega_{n+1}^2-4\omega_n\omega_{n+1}+12\omega_n^2,\\
&c_n=2A^2+10B^2=2\omega_{n+1}^2+10\omega_n^2\\
&d_n= A^2-9AB-B^2=\omega_{n+1}^2-9\omega_n\omega_{n+1}-9\omega_n^2.
\end{split}
\end{align}
satisfying
\begin{align*}
a_n^3+b_n^3=c_n^3+d_n^3.
\end{align*}
Now we solve the recurrence relation of $\omega_n$ to see
\begin{align*}
&\omega_n=\frac{1}{\sqrt{5}}\bigg(\bigg(\frac{1+\sqrt{5}}{2}\bigg)^n+\bigg(\frac{1-\sqrt{5}}{2}\bigg)^n\bigg),\\
&\omega_n^2=\frac{1}{5}\bigg(\bigg(\frac{3+\sqrt{5}}{2}\bigg)^n+\bigg(\frac{3-\sqrt{5}}{2}\bigg)^n-2(-1)^n\bigg),\\&\omega_n\omega_{n+1}=\frac{1}{5}\bigg(\bigg(\frac{1+\sqrt{5}}{2}\bigg)^{2n+1}+\bigg(\frac{1-\sqrt{5}}{2}\bigg)^{2n+1}-(-1)^n\bigg).\\
\end{align*}
This leads us to obtain the following generating functions:
\begin{align}\label{generating function of omega_n}
\begin{split}
&\sum_{n=0}^\infty \omega_n^2x^n=\frac{x-x^2}{1-2x-2x^2+x^3},\\
&\sum_{n=0}^\infty \omega_n\omega_{n+1}x^n=\frac{x}{1-2x-2x^2+x^3}.
\end{split}
\end{align}
We now make use of \eqref{generating function of omega_n} in \eqref{a_n, b_n, c_n, d_n} to get the generating functions for $a_n, b_n, c_n$ and $d_n.$ This completes the proof.
\end{proof}
\begin{proof}[Theorem $\ref{gen -9^3n}$][]
We define a sequence $\{\omega_n\}$ satisfying the following recurrence relation:
\begin{align*}
\omega_{n+2}=9\omega_n-7\omega_{n+1},\hspace{3mm} \textrm{with}\hspace{3mm} \omega_0=0,\hspace{1mm} \omega_1=1.
\end{align*}
Then one can verify that $\omega_n$ satisfies 
\begin{align}
\omega_{n+1}^2-\omega_n\omega_{n+2} = (-9)^n.
\end{align}
Now, to synchronize the sequence $\omega_n$ with \eqref{after_variable_change}, we substitute
\begin{align}
 A=\omega_{n+1},~~B=\omega_n.
\end{align}
This leads us to obtain
\begin{align}
A^2+7AB-9B^2 &= \omega_{n+1}^2 + 7\omega_n\omega_{n+1} - 9\omega_n^2 \nonumber\\ &=  \omega_{n+1}^2 - \omega_n(9\omega_n - 7\omega_{n+1})\nonumber \\&=  \omega_{n+1}^2 - \omega_n \omega_{n+2} = (-9)^n.
\end{align}
We write the remaining terms of \eqref{after_variable_change} as
\begin{align}\label{a_n, b_n, c_n for 9^n}
\begin{split}
&a_n=2A^2+10B^2=2\omega_{n+1}^2+10\omega_n^2,\\
&b_n=A^2-9AB-B^2=\omega_{n+1}^2-9\omega_n\omega_{n+1}-\omega_n^2, \\
&c_n=2A^2-4AB+12B^2=2\omega_{n+1}^2-4\omega_n\omega_{n+1}+12\omega_n^2,
\end{split}
\end{align}
and they satisfy 
\begin{align*}
a_n^3+b_n^3=c_n^3+((-9)^n)^3.
\end{align*} Solving recurrence relation for $\omega_n$, we have  
\begin{align*}
&\omega_n=\frac{1}{\sqrt{85}}\bigg(\frac{1}{2^n}(\sqrt{85}-7)^n-\frac{1}{2^n}(-\sqrt{85}-7)^n\bigg),\\
&\omega_n^2=\frac{1}{85}\bigg(\frac{1}{2^n}(7\sqrt{85}+67)^n+\frac{1}{2^n}(-7\sqrt{85}+67)^n-2(-9)^n\bigg),\\
&\omega_n\omega_{n+1}=\frac{1}{85}\bigg(\frac{1}{2^{2n+1}}(\sqrt{85}-7)^{2n+1}+\frac{1}{2^{2n+1}}(-\sqrt{85}-7)^{2n+1}+7(-9)^n\bigg).
\end{align*} This leads us to obtain the following generating functions:
\begin{align}\label{generating function for w_n3}
\begin{split}
&\sum_{n=0}^\infty \omega_n^2x^n=\frac{x-9x^2}{1-58x-522x^2+729x^3},\\
&\sum_{n=0}^\infty \omega_n\omega_{n+1}x^n=\frac{-7x}{1-58x-522x^2+729x^3}.
\end{split}
\end{align}
Utilizing \eqref{generating function for w_n3} in \eqref{a_n, b_n, c_n for 9^n} completes the proof of first equality of \eqref{abc}. Now expanding the left hand sides of \eqref{generating functions 9^n} around $x=\infty$ and collecting the coefficient of $x^{-n-1}$, we get
\begin{align*}
\alpha_n=\frac{-162}{85\sqrt{85}}\big[(71\eta-11)\eta^n - (71\xi-11)\xi^n\big]+\frac{16}{765}\left(\frac{-1}{9}\right)^n, \\
\beta_n=\frac{648}{595\sqrt{85}}\big[(59\eta+16)\eta^n - (59\xi+16)\xi^n\big]-\frac{43}{765}\left(\frac{-1}{9}\right)^n, \\
\gamma_n=\frac{-486}{595\sqrt{85}}\big[(146\eta-31)\eta^n - (146\xi-31)\xi^n\big]-\frac{16}{765}\left(\frac{-1}{9}\right)^n,
\end{align*}
where $\eta$ and $\xi$ are the roots of the polynomial equation $81x^2-67x+1=0$. Replacing $n$ by $n-1$, we get the required result.
\end{proof}
\begin{proof}[Theorem $\ref{gen 6^3n}$][]
The proof of this theorem goes along the same line as in Theorem \ref{gen -9^3n}. Here the terms of the sequence $\{\omega_n\}$ can be derived from the following recurrence relation:
\begin{align*}
\omega_{n+2}=-6\omega_n+2\omega_{n+1},\hspace{3mm} \omega_0=0,\hspace{1mm} \omega_1=1,
\end{align*}
and satisfies
\begin{align*}
2\omega_{n+1}^2-2\omega_n\omega_{n+2}=2\cdot6^n.
\end{align*}
Here we utilize \eqref{after_variable_change} by substituting $A=\omega_{n+1},~~B=\omega_n,$ where 
\begin{align}\label{a_n, b_n, c_n for 2.6^n}
\begin{split}
&a_n=2A^2+10B^2=2\omega_{n+1}^2+10\omega_n^2,\\
&b_n=A^2-9AB-B^2=\omega_{n+1}^2-9\omega_n\omega_{n+1}-9\omega_n^2,\\
&c_n=A^2+7AB-9B^2=\omega_{n+1}^2+7\omega_n\omega_{n+1}-9\omega_n^2,
\end{split}
\end{align}
and
\begin{align*}
2A^2-4AB+12B^2=2(6)^n.
\end{align*}
Then we have
\begin{align*}
a_n^3+b_n^3=c_n^3+2^3(6)^{3n}.
\end{align*}
We solve the recurrence relation of $\omega_n$ to get,
\begin{align*}
&\omega_n=\frac{i}{2 \sqrt{5}}\bigg((1-i\sqrt{5})^n-(1+i\sqrt{5})^n\bigg),\\
&\omega_n^2=\frac{-1}{20}\bigg((-4-2i\sqrt{5})^n+(-4+2i\sqrt{5})^n -2(6)^n\bigg),\\
&\omega_n\omega_{n+1}=\frac{-1}{20}\bigg(1-i\sqrt{5})^{2n+1}+(1+i\sqrt{5})^{2n+1}-2(6)^n\bigg).
\end{align*}
Then the generating functions for $\omega_n^2$ and $\omega_n \omega_{n+1}$ are
\begin{align*}
&\sum_{n=0}^\infty \omega_n^2x^n=\frac{x+6x^2}{1+2x-12x^2-216x^3},\\
&\sum_{n=0}^\infty \omega_n\omega_{n+1}x^n=\frac{2x}{1+2x-12x^2-216x^3}.
\end{align*}
Utilizing these generating functions in \eqref{a_n, b_n, c_n for 2.6^n}, we get the required result.
\end{proof}


\begin{proof}[Theorem \ref{Fibonacci_power_4_s_t}][]
The coefficients of the power series obtained in \eqref{fibonacci a_n...f_n} gives solution to the Diophantine equation $A^4 +B^4+C^4+D^4+E^4=F^4$. To prove this theorem, we make use of the following parametric solution given by Ramanujan \cite[p.~386]{ramanujantifr}. If $s$ and $t$ are arbitrary, then
\begin{align}\label{power_4_s_t}
\begin{split}
(8s^2 + 40st &-24t^2)^4 + (6s^2 - 44st -18t^2)^4 + (14s^2 -4st -42t^2)^4 \\ &+ (9s^2 + 27t^2)^4 + (4s^2 + 12t^2)^4 = (15s^2 + 45t^2)^4.
\end{split}
\end{align}
The proof goes in the same vein as done in previous theorems. Here we take our sequence $\{\omega_n \}$ to be the Fibonacci sequence defined as in \eqref{Fibonnaci_sequence}. Now to synchronize this with \eqref{power_4_s_t}, we substitute,
\begin{align*}
s=\omega_{n+1}\hspace{3mm} \textrm{and} \hspace{3mm}t=\omega_n.
\end{align*}
Then we represent $a_n, b_n, c_n, d_n, e_n$ and $f_n$ as
\begin{align}\label{a_n,b_n...f_n,fibonacci}
\begin{split}
&a_n= 8s^2 + 40st -24t^2 = 8\omega_{n+1}^2+40\omega_n\omega_{n+1}-24\omega_n^2,   \\
&b_n= 6s^2 - 44st -18t^2 =  6\omega_{n+1}^2-44\omega_n\omega_{n+1}-18\omega_n^2,   \\
&c_n= 14s^2 -4st -42t^2 =  14\omega_{n+1}^2-4\omega_n\omega_{n+1}-42\omega_n^2,   \\
&d_n= 9s^2 + 27t^2= 9\omega_{n+1}^2+27\omega_n^2,   \\
&e_n= 4s^2 + 12t^2= 4\omega_{n+1}^2+12\omega_n^2,   \\
&f_n= 15s^2 + 45t^2 = 15\omega_{n+1}^2+45\omega_n^2,
\end{split}
\end{align}
satisfying
\begin{align*}
a_n^4 + b_n^4 + c_n^4 + d_n^4 + e_n^4 = f_n^4. 
\end{align*}
The generating functions for $\omega_n^2$ and $\omega_n \omega_{n+1}$ are given in \eqref{generating function of omega_n}. Utilizing \eqref{generating function of omega_n} in \eqref{a_n,b_n...f_n,fibonacci},  one can complete the proof.
\end{proof}
\begin{proof}[Theorem \ref{3^n.8_power_4}][]
Let us consider a sequence $\{\omega_n\}$ defined as 
\begin{align}\label{w_n 2.8}
\omega_{n+2}=3\omega_n-5\omega_{n+1},\hspace{3mm} \omega_0=0,\hspace{1mm} \omega_1=1.
\end{align}
This sequence satisfies a recurrence relation $$g_n = (-3)^n g_0 = 8(-3)^n ,$$
where we define $g_n$ as $$g_n = 8\omega_{n+1}^2 - 8\omega_n \omega_{n+2}.$$
Now, in accordance with \eqref{power_4_s_t}, we substitute 
\begin{align*}
s=\omega_{n+1} \hspace{5mm} \textrm{and} \hspace{5mm} t=\omega_n,
\end{align*}
and we write
\begin{align}\label{a_n,b_n...f_n,8.3^n}
\begin{split}
&a_n= 6s^2 - 44st -18t^2 =  6\omega_{n+1}^2-44\omega_n\omega_{n+1}-18\omega_n^2,   \\
&b_n= 14s^2 -4st -42t^2 =  14\omega_{n+1}^2-4\omega_n\omega_{n+1}-42\omega_n^2,   \\
&c_n= 9s^2 + 27t^2= 9\omega_{n+1}^2+27\omega_n^2,   \\
&d_n= 4s^2 + 12t^2= 4\omega_{n+1}^2+12\omega_n^2,   \\
&e_n= 15s^2 + 45t^2 = 15\omega_{n+1}^2+45\omega_n^2,
\end{split}
\end{align}
satisfying
\begin{align*}
a_n^4 + b_n^4 + c_n^4 + d_n^4 = e_n^4 - (8\times (-3)^n)^4.
\end{align*}
On solving \eqref{w_n 2.8}, we get the generating functions for $w_n$ to be
\begin{align*}
&\sum_{n=0}^\infty \omega_n^2x^n=\frac{x-3x^2}{1-28x-84x^2+27x^3},\\
&\sum_{n=0}^\infty \omega_n\omega_{n+1}x^n=\frac{-5x}{1-28x-84x^2+27x^3}.
\end{align*}
Make use of these generating functions in \eqref{a_n,b_n...f_n,8.3^n} to get the desired result.
\end{proof}
\begin{proof}[Theorem \ref{fibonacci_power_4_m_n}][]
The proof of this theorem goes along the lines of Theorem \ref{Fibonacci_power_4_s_t}. The key difference here is that we use a different parametric solution of Diophantine equation $A^4 +B^4+C^4+D^4+E^4=F^4$ given by Ramanujan \cite[p.~386]{ramanujantifr}. Basically, for arbitrary $m$ and $n$, we have
\begin{align}\label{power_4_m_n}
\begin{split}
(4m^2-&12n^2)^4 + (3m^2 + 9n^2)^4 + (2m^2 - 12mn -6n^2)^4 \\& + (4m^2+12n^2)^4 + (2m^2 + 12mn -6n^2)^4 = (5m^2 + 15n^2)^4.
\end{split}
\end{align}
Here again, we consider our sequence to be Fibonacci sequence. As the proof goes along the lines of Theorem \ref{fibonacci},  so we left it for the readers. 
\end{proof}

\begin{proof}[Theorem \ref{2.3^n_Power_4}][]
Here we consider our sequence $\{\omega_n\}$ to be
\begin{align*}
\omega_{n+2}=6\omega_{n+1}+3\omega_{n},\hspace{3mm} \omega_0=0,\hspace{1mm} \omega_1=1.
\end{align*}
Now let $$g_n=2\omega_{n+1}^2-2\omega_n \omega_{n+2},$$ then it satisfies $$g_n=(-3)^ng_0=2\times(-3)^n.$$ The rest of the proof goes along the same lines as in Theorem \ref{3^n.8_power_4},  so we omit the proof.
\end{proof}

{\bf Acknowledgement:}
We sincerely thank the anonymous referee for carefully reading our manuscript and for giving valuable suggestions.
The authors would like to thank Dr. Ajai Choudhry, Dr. Pramod Eyyunni and Dr.  Bibekananda Maji for their valuable suggestions throughout this project.  The first author's research is supported by the CSIR-UGC Fellowship, Govt. of India, whereas the second author's research is given by the PMRF Fellowship, Govt. of India, grant number 2101705. We sincerely thank our respective agencies for their generous support. We would also like to thank Bhaskaracharya Mathematics Laboratory and Brahmagupta Mathematics Library of the Department of Mathematics, IIT Indore, supported by DST FIST Project (File No.: SR/FST/MS-I/2018/26).


\begin{thebibliography}{99}

\bibitem{LN4} G.~E.~Andrews and B.~C.~Berndt, \emph{Ramanujan’s Lost Notebook, Part IV}, Springer, New York, 2013.

\bibitem{RN4} B.~C.~Berndt, \emph{Ramanujan’s Notebook, Part IV}, Springer, New York, 1994.

\bibitem{chen} K.~W.~Chen, \emph{Extension of an amazing identity of Ramanujan,} Fibonacci Quart. \textbf{50} (2012), 227--230.

\bibitem{DGR} J.~H.~Dinitz, R.~Games and R.~Roth, \emph{Seeds for generalized taxicab numbers,} J. Integer Seq. \textbf{22} (2019), Article 19.3.3.

\bibitem{hir 1} J.~H.~Han and M.~D.~Hirschhorn, \emph{Another look at an amazing identity of Ramanujan}, Math. Magazine \textbf{79} (2006), 302--304.

\bibitem{hir 2} M.~D.~Hirschhorn, \emph{An amazing identity of Ramanujan}, Math. Magazine \textbf{68} (1995), 199--201.

\bibitem{hir 3} M.~D.~Hirschhorn, \emph{A proof in the spirit of Zeilberger of an amazing identity of Ramanujan}, Math. Magazine \textbf{69} (1996), 267--269.

\bibitem{hir 4} M.~D.~Hirschhorn, \emph{Ramanujan and Fermat's last theorem}, Austral. Math. Soc. Gaz. \textbf{31} (2004), 256--257.

\bibitem{JMc} J.~McLaughlin, \emph{An identity motivated by an amazing identity of Ramanujan,} The Fibinacci Quaterly \textbf{48.1} (2010), 34--38.

\bibitem{VN} V.~Ni\c{t}ic\u{a}, \emph{About some relatives of the taxicab number}, J. Integer seq. \textbf{21} (2018), Article 18.9.4.

\bibitem{ramanujantifr} S.~Ramanujan, \emph{Notebooks of Srinivasa Ramanujan, Vol. II}, Tata Institute of Fundamental Research, Mumbai, 2012.

\bibitem{ramanujanquestion} S.~Ramanujan \emph{Question} 441, J. Indian Math. Soc. \textbf{5} (1913), 39; solution in \textbf{6} (1914), 226--227.

\bibitem{lost notebook} S.~Ramanujan, \emph{The Lost Notebook and Other Unpublished Papers,} New Delhi, Narosa, 1988.


\bibitem{silverman} J.~H.~Silverman, \emph{Taxicabs and sums of two cubes,} Amer. Math. Monthly, \textbf{100} (1993), 331--340.

\end{thebibliography}
\end{document}